# Stability of a Nonlinear Axially Moving String With the Kelvin-Voigt Damping *

S. M. Shahruz

Berkeley Engineering Research Institute
P. O. Box 9984
Berkeley, California 94709

## Abstract

In this paper, a nonlinear axially moving string with the Kelvin-Voigt damping is considered. It is proved that the string is stable, i.e., its transversal displacement converges to zero when the axial speed of the string is less than a certain critical value. The proof is established by showing that a Lyapunov function corresponding to the string decays to zero exponentially. It is also shown that the string displacement is bounded when a bounded distributed force is applied to it transversally. Furthermore, a few open problems regarding the stability and stabilization of strings with the Kelvin-Voigt damping are stated.

## 1. Introduction

In the past few decades, the nonlinear dynamics and control of axially moving strings have been studied by researchers; see, e.g., Refs. [1-13]. Moving string-like continua, such as threads, wires, tethers, cables, magnetic tapes, belts, band-saws, chains, etc., are used in many applications.

---

* This paper has been submitted to a journal for publication.

In most applications, a string-like continuum is pulled at a constant speed through two fixed eyelets, as shown in Fig. 1. The distance between the eyelets in this figure is assumed to be $1$. There are two commonly used models that represent the dynamics of moving strings. The first model, known as the moving Kirchhoff string, is represented by

$$y_{tt}(x, t) + 2v\, y_{xt}(x, t) = \left[ 1 - v^2 + b \int_0^1 y_x^2(x, t)\, dx \right] y_{xx}(x, t)\,, \tag{1}$$

for all $x \in (0, 1)$ and $t \geq 0$; see, e.g., Refs. [5, 7, 8, 12]. In Eq. (1), $y(\cdot, \cdot) \in \mathbb{R}$ denotes the transversal displacement of the string in the $Y$-direction, $y_x := \partial y / \partial x$, $y_{xx} := \partial^2 y / \partial x^2$, $y_{tt} := \partial^2 y / \partial t^2$, and $y_{xt} := \partial^2 y / \partial x \partial t$; $b > 0$ is a constant real number, and $v \geq 0$ is proportional to the speed of the string through the eyelets. For realistic situations, $0 \leq v < 1$.

The second model is represented by

$$y_{tt}(x, t) + 2v\, y_{xt}(x, t) = \left[ 1 - v^2 + \frac{3}{2}\, b\, y_x^2(x, t) \right] y_{xx}(x, t)\,, \tag{2}$$

for all $x \in (0, 1)$ and $t \geq 0$; see, e.g., Refs. [1, 6, 10]. In Eq. (2), $y(\cdot, \cdot) \in \mathbb{R}$ is the string displacement and its derivatives with respect to $x$ and $t$ are readily identified; $b > 0$ is a constant real number, and $0 \leq v < 1$ is proportional to the string speed.

The boundary conditions of the strings represented by Eqs. (1) and (2) can be specified depending on whether the eyelets are both fixed or a boundary control is applied at one eyelet. The initial displacement and velocity of the strings can be specified so that the initial boundary value problems corresponding to the strings would be well defined.

There is yet another nonlinear model of moving strings obtained by using the Kelvin-Voigt constitutive law. This model is given by

$$\begin{aligned} y_{tt}(x, t) + 2\delta\, y_t(x, t) + 2v\, y_{xt}(x, t) &= \left[ 1 - v^2 + \frac{3}{2}\, b\, y_x^2(x, t) \right] y_{xx}(x, t) \\ &\quad + 2\eta\, y_x(x, t)\, y_{xx}(x, t)\, y_{xt}(x, t) + \eta\, y_x^2(x, t)\, y_{xxt}(x, t)\,, \end{aligned} \tag{3}$$

for all $x \in (0, 1)$ and $t \geq 0$; see, e.g., Refs. [9, 13]. In Eq. (3), $y(\cdot, \cdot) \in \mathbb{R}$ is the string displacement, $y_t := \partial y / \partial t$, and other derivatives of $y(\cdot, \cdot)$ with respect to $x$ and $t$ are readily identified; $b > 0$, $\delta > 0$ is the viscous damping coefficient, $\eta > 0$ is known as the Kelvin-Voigt damping, and $0 \leq v < 1$ is proportional to the string speed. (It is remarked that $\eta$ may not be called the Kelvin-Voigt damping by some authors. However, since it corresponds to the Kelvin-Voigt constitutive law and is related to energy dissipation in the string, it will be called as such in this paper.)

For fixed eyelets, the boundary conditions corresponding to the string represented by Eq. (3) are

$$y(0, t) = y(1, t) = 0 , \tag{4}$$

for all $t \geq 0$ . The initial displacement and velocity of the string are, respectively,

$$y(x, 0) = f(x) , \quad y_t(x, 0) = g(x) , \tag{5}$$

for all $x \in (0, 1)$ . It is assumed that $f \in C^1[0, 1]$ and that at least one of the functions $f$ or $g$ is not identically zero over $[0, 1]$ .

In this paper, the goals are: (i) to show that the nonlinear axially moving string represented by Eqs. (3), (4), and (5) is stable; that is, $y(x, t) \to 0$ as $t \to \infty$ for all $x \in (0, 1)$ ; (ii) to establish the bounded-input bounded-output (BIBO) stability of the string; that is, to show that the string displacement is bounded when a bounded distributed force is applied to it transversally. Extensive literature survey proved that these stability results are not available. Furthermore, a few open problems regarding the stability and stabilization of strings with the Kelvin-Voigt damping are stated.

## 2. Stability of the Moving String

The plan to establish the stability of the nonlinear moving string represented by Eqs. (3), (4), and (5) is as follows. A Lyapunov function is proposed for the string. This function is denoted by $t \mapsto V(t)$ . It will be shown that $V(\cdot)$ tends to zero exponentially, thereby the stability of the string will be established. It is remarked that the Lyapunov function $V(\cdot)$ was obtained after so much effort.

The scalar-valued function $V(\cdot)$ is defined as

$$V(t) = E(t) + \delta \int_0^1 [y(x, t)\, y_t(x, t) + \delta\, y^2(x, t) + \frac{\eta}{4}\, y_x^4(x, t)]\, dx , \tag{6}$$

for all $t \geq 0$ , where

$$E(t) = \int_0^1 \frac{1}{2}\, y_t^2(x,\, t)\; dx + \int_0^1 \frac{1}{2}\,(1 - \nu^2)\; y_x^2(x,\, t)\; dx + \int_0^1 \frac{b}{8}\; y_x^4(x,\, t)\; dx\,. \tag{7}$$

From Eqs. (6) and (7), it follows that

$$V(t) = \int_0^1 \frac{1}{2}\,[y_t(x,\, t) + \delta\; y(x,\, t)]^2 dx + \int_0^1 \frac{\delta^2}{2}\; y^2(x,\, t)\; dx$$

$$+ \int_0^1 \frac{1}{2}\,(1 - \nu^2)\; y_x^2(x,\, t)\; dx + \int_0^1 (\frac{b}{8} + \frac{\delta\eta}{4})\; y_x^4(x,\, t)\; dx\,, \tag{8}$$

for all $t \geq 0$ . From Eqs. (8) and (5), it is concluded that

$$V(0) = \int_0^1 \frac{1}{2}\,[g(x) + \delta\;\, f(x)]^2 dx + \int_0^1 \frac{\delta^2}{2}\;\, f^2(x)\; dx$$

$$+ \int_0^1 \frac{1}{2}\,(1 - \nu^2)\;\, f_x^2(x)\; dx + \int_0^1 (\frac{b}{8} + \frac{\delta\eta}{4})\;\, f_x^4(x)\; dx\,, \tag{9}$$

where $f_x(x) := df(x)/dx$ . Since at least one of the functions $f$ or $g$ in Eq. (5) is not identically zero over $[0, 1]$ , it is evident that $V(0) > 0$ .

A property of $V(\cdot)$ is now proved.

**Lemma 2.1:** The function $V(\cdot)$ satisfies

$$0 \leq V(t) \leq K\;\, E(t)\,, \tag{10}$$

for all $t \geq 0$ , where

$$K := 1 + \delta\; \max\left\{ \frac{1 + 2\delta/\pi}{\pi(1 - \nu^2)},\, \frac{2\eta}{b} \right\}. \tag{11}$$

**Proof:** From Eq. (8), it is clear that $V(t) \geq 0$ for all $t \geq 0$ . For the first two integral terms in Eq. (6) which are multiplied by $\delta$ , the following inequalities hold, respectively:

$$\int_0^1 y\, y_t\, dx \le \frac{\pi}{2} \int_0^1 y^2\, dx + \frac{1}{2\pi} \int_0^1 y_t^2\, dx$$

$$\le \frac{1}{2\pi} \int_0^1 y_x^2\, dx + \frac{1}{2\pi} \int_0^1 y_t^2\, dx = \frac{1}{\pi(1-v^2)} \int_0^1 \frac{1}{2} (1-v^2)\, y_x^2\, dx + \frac{1}{\pi} \int_0^1 \frac{1}{2}\, y_t^2\, dx\,, \quad \text{(12a)}$$

$$\int_0^1 \delta\, y^2\, dx \le \frac{\delta}{\pi^2} \int_0^1 y_x^2\, dx = \frac{2\delta/\pi}{\pi(1-v^2)} \int_0^1 \frac{1}{2} (1-v^2)\, y_x^2\, dx\,, \quad \text{(12b)}$$

for all $t \ge 0$, where the arguments $(x, t)$ of functions are deleted. In deriving Eq. (12), Scheeffer's inequality, which is a Poincare-type inequality, is used; see, e.g., Ref. [14, p. 67]. Using Eq. (12) and the fact that $v < 1$, it follows that

$$\int_0^1 (y\, y_t + \delta\, y^2 + \frac{\eta}{4}\, y_x^4)\, dx \le \max\left\{ \frac{1}{\pi}, \frac{1+2\delta/\pi}{\pi(1-v^2)}, \frac{2\eta}{b} \right\}$$

$$\times \left( \int_0^1 \frac{1}{2}\, y_t^2\, dx + \int_0^1 \frac{1}{2} (1-v^2)\, y_x^2\, dx + \int_0^1 \frac{b}{8}\, y_x^4\, dx \right)$$

$$\le \max\left\{ \frac{1+2\delta/\pi}{\pi(1-v^2)}, \frac{2\eta}{b} \right\} E(t)\,, \quad \text{(13)}$$

for all $t \ge 0$, where $E(\cdot)$ is that in Eq. (7). Using Eq. (13) in Eq. (6), it is concluded that Eq. (10) holds. □

Next, several identities and an inequality for functions satisfying Eq. (4) are proved.

**Lemma 2.2:** Let $y(\cdot, \cdot)$ satisfy the boundary conditions in Eq. (4). Then,

$$\int_0^1 2\, y_t\, y_{xt}\, dx = 0\,, \quad \text{(14a)}$$

$$\int_0^1 (y_{xx}\, y_t + y_{xt}\, y_x)\, dx = 0\,, \quad \text{(14b)}$$

$$\int_0^1 (3\, y_x^2\, y_{xx}\, y_t + y_x^3\, y_{xt})\, dx = 0\,, \quad \text{(14c)}$$

$$\int_0^1 y_t\,(y_x^2\ y_{xt})_x\ dx = -\int_0^1 y_x^2\ y_{xt}^2\ dx\ , \tag{14d}$$

$$-\,2\int_0^1 y\ y_{xt}\ dx \le \int_0^1 y_x^2\ dx + \int_0^1 y_t^2\ dx\ , \tag{14e}$$

$$\int_0^1 y\ y_{xx}\ dx = -\int_0^1 y_x^2\ dx\ , \tag{14f}$$

$$\int_0^1 3\ y\ y_x^2\ y_{xx}\ dx = -\int_0^1 y_x^4\ dx\ , \tag{14g}$$

$$\int_0^1 y\,(y_x^2\ y_{xt})_x\ dx = -\,\frac{1}{4}\int_0^1 (y_x^4)_t\ dx\ , \tag{14h}$$

for all $t \ge 0$ .

**Proof:** From Eq. (4), it follows that $y_t(0, t) = y_t(1, t) = 0$ for all $t \ge 0$ . Hence,

$$\int_0^1 2\ y_t\ y_{xt}\ dx = \int_0^1 (y_t^2)_x\ dx = 0\ , \tag{15}$$

for all $t \ge 0$ . Thus, Eq. (14a) holds.

Using $y_t(0, t) = y_t(1, t) = 0$ for all $t \ge 0$ , it is concluded that

$$\int_0^1 (y_{xx}\ y_t + y_{xt}\ y_x)\ dx = \int_0^1 (y_x\ y_t)_x\ dx = 0\ , \tag{16a}$$

$$\int_0^1 (3\ y_x^2\ y_{xx}\ y_t + y_x^3\ y_{xt})\ dx = \int_0^1 (y_x^3\ y_t)_x\ dx = 0\ , \tag{16b}$$

for all $t \ge 0$ . Thus, Eqs. (14b) and (14c) hold.

Next,

$$\int_0^1 y_t\,(y_x^2\ y_{xt})_x\ dx = -\int_0^1 y_x^2\ y_{xt}^2\ dx\ , \tag{17}$$

for all $t \geq 0$, where the last identity is obtained by integrating by parts and using $y_t(0, t) = y_t(1, t) = 0$ for all $t \geq 0$. Thus, Eq. (14d) holds.

The proof of Eq. (14e) is as follows:

$$-\int_0^1 2\, y\, y_{xt}\, dx = \int_0^1 2\, y_x\, y_t\, dx \leq \int_0^1 y_x^2\, dx + \int_0^1 y_t^2\, dx\,, \tag{18}$$

for all $t \geq 0$, where the first identity is obtained by integrating by parts and using the boundary conditions in Eq. (4).

The proof of Eq. (14f) is established by integrating by parts and using the boundary conditions in Eq. (4); so are established those of Eqs. (14g) and (14h). The proofs of these two last identities are, respectively, given in the following:

$$\int_0^1 3\, y\, y_x^2\, y_{xx}\, dx = \int_0^1 y\, (y_x^3)_x\, dx = -\int_0^1 y_x^4\, dx\,, \tag{19a}$$

$$\int_0^1 y\, (y_x^2\, y_{xt})_x\, dx = -\int_0^1 y_x^3\, y_{xt}\, dx = -\frac{1}{4}\int_0^1 (y_x^4)_t\, dx\,, \tag{19b}$$

for all $t \geq 0$. □

Using these preliminary results, it is next proved that $V(\cdot)$ tends to zero exponentially.

**Theorem 2.3:** If

$$0 \leq \nu < \nu_c := \frac{\sqrt{5} - 1}{2}\,, \tag{20}$$

then the function $V(\cdot)$ in Eq. (6) along the solution of Eqs. (3), (4), and (5) satisfies

$$V(t) \leq V(0)\, \exp\left( -\frac{2\delta(1 - \nu - \nu^2)}{K(1 - \nu^2)}\, t \right), \tag{21}$$

for all $t \geq 0$, where $K$ is that in Eq. (11). That is, $V(\cdot)$ tends to zero exponentially.

**Proof:** From Eq. (8), it follows that

$$\dot{V}(t) = \int_0^1 (y_t + \delta\, y)\,(y_{tt} + \delta\, y_t)\, dx + \int_0^1 \delta^2\, y\, y_t\, dx$$

$$+ \int_0^1 (1 - v^2)\, y_x\, y_{xt}\, dx + \int_0^1 \frac{b}{2}\, y_x^3\, y_{xt}\, dx + \int_0^1 \frac{\delta\eta}{4}\, (y_x^4)_t\, dx\,, \tag{22}$$

for all $t \geq 0$. Substituting $y_{tt}$ from Eq. (3) into Eq. (22) and noting that $2\, y_x\, y_{xx}\, y_{xt} + y_x^2\, y_{xxt} = (y_x^2\, y_{xt})_x$, it is concluded that

$$\dot{V}(t) = -\,\delta \int_0^1 y_t^2\, dx - v \int_0^1 2\, y_t\, y_{xt}\, dx$$

$$+ (1 - v^2) \int_0^1 (y_{xx}\, y_t + y_{xt}\, y_x)\, dx + \frac{b}{2} \int_0^1 (3\, y_x^2\, y_{xx}\, y_t + y_x^3\, y_{xt})\, dx$$

$$+ \eta \int_0^1 y_t\, (y_x^2\, y_{xt})_x\, dx - 2\delta v \int_0^1 y\, y_{xt}\, dx + \delta(1 - v^2) \int_0^1 y\, y_{xx}\, dx$$

$$+ \frac{b\delta}{2} \int_0^1 3\, y\, y_x^2\, y_{xx}\, dx + \delta\eta \int_0^1 y\, (y_x^2\, y_{xt})_x\, dx + \frac{\delta\eta}{4} \int_0^1 (y_x^4)_t\, dx\,, \tag{23}$$

for all $t \geq 0$. Using Eqs. (14a)-(14h) in Eq. (23), the following inequality is obtained:

$$\dot{V}(t) \leq -\,\delta(1 - v) \int_0^1 y_t^2\, dx - \delta(1 - v - v^2) \int_0^1 y_x^2\, dx - \frac{b\delta}{2} \int_0^1 y_x^4\, dx - \eta \int_0^1 y_x^2\, y_{xt}^2\, dx\,, \tag{24}$$

for all $t \geq 0$. By neglecting the last term on the right-hand side of Eq. (24), rearranging the other terms, and using Eq. (7), it follows that

$$\dot{V}(t) \leq -\,2\delta \left( \frac{1 - v - v^2}{1 - v^2} \right) E(t) - \delta \left( \frac{v^3}{1 - v^2} \int_0^1 y_t^2\, dx + \frac{b}{4} \left( \frac{1 + v - v^2}{1 - v^2} \right) \int_0^1 y_x^4\, dx \right), \tag{25}$$

for all $t \geq 0$, where the coefficients of the last two integral terms are non-positive for $0 \leq v < 1$. Therefore,

$$\dot{V}(t) \leq -\,2\delta \left( \frac{1 - v - v^2}{1 - v^2} \right) E(t)\,, \tag{26}$$

for all $t \geq 0$, where by Eq. (20), the coefficient of $E(\cdot)$ is negative. Using Eq. (10) in Eq. (26), the following differential inequality is obtained:

$$\dot{V}(t) \leq -\frac{2\delta(1 - v - v^2)}{K(1 - v^2)} V(t) , \tag{27}$$

for all $t \geq 0$, with the initial condition $V(0) > 0$ given in Eq. (9). By a comparison theorem in Ref. [15, p. 3] or Ref. [16, p. 2], it is concluded that Eq. (21) holds. □

Finally, the stability of the string can be established.

**Corollary 2.4:** The solution of Eqs. (3), (4), and (5), $y(x, t) \rightarrow 0$ as $t \rightarrow \infty$ for all $x \in [0, 1]$.

**Proof:** For the string represented by Eqs. (3), (4), and (5), the Lyapunov function in Eq. (6) is chosen. By Theorem 2.3, the function $V(\cdot)$ tends to zero exponentially. From Eq. (8), it follows that $y_x(x, t) \rightarrow 0$ as $t \rightarrow \infty$ for all $x \in [0, 1]$. Since $y(0, t) = 0$ for all $t \geq 0$, it is concluded that $y(x, t) \rightarrow 0$ as $t \rightarrow \infty$ for all $x \in [0, 1]$. □

In establishing the stability of the nonlinear moving string, the viscous damping in the model played an important role; so did the the Lyapunov function, which was obtained after so many trials. By Eq. (20), the stability of the string is guaranteed when the string speed is less than the critical speed $v_c \approx 0.6$.

## 3. Bounded-Input Bounded-Output Stability

In this section, the string represented by Eqs. (3) and (4) is reconsidered. It is assumed that the initial displacement and velocity of the string are zero, however, an external distributed force is applied to it transversally. In this case, the string is represented by

$$y_{tt}(x, t) + 2\delta\, y_t(x, t) + 2v\, y_{xt}(x, t) = \left[ 1 - v^2 + \frac{3}{2} b\, y_x^2(x, t) \right] y_{xx}(x, t)$$

$$+ 2\eta\, y_x(x, t)\, y_{xx}(x, t)\, y_{xt}(x, t) + \eta\, y_x^2(x, t)\, y_{xxt}(x, t) + F(x, t) , \tag{28}$$

for all $x \in (0, 1)$ and $t \geq 0$, where $F(x, t) \in \mathbb{R}$ is the applied force, the boundary conditions are the same as those in Eq. (4), and the initial conditions are

$$y(x,0) = 0\,, \qquad y_t(x,0) = 0\,. \tag{29}$$

In this section, the goal is to show that the string displacement is bounded when the external force is bounded. To define boundedness precisely, two function spaces are introduced:

**(i)** Let $X_2$ denote the space of functions $u : (0,1) \times \mathbb{R}_+ \to \mathbb{R}$, given by $u(x,t)$, for which $\|u\|_{X_2} := \sup_{t \geq 0} \left[ \int_0^1 |u(x,t)|^2 \, dx \right]^{1/2} < \infty$. A $u \in X_2$ is said to be $X_2$-bounded.

**(ii)** Let $X_\infty$ denote the space of functions $u : (0,1) \times \mathbb{R}_+ \to \mathbb{R}$, given by $u(x,t)$, for which $\|u\|_{X_\infty} := \sup_{t \geq 0} \sup_{x \in (0,1)} |u(x,t)| < \infty$. A $u \in X_\infty$ is said to be $X_\infty$-bounded.

It is clear that $X_\infty \subset X_2$ since $L_\infty(0,1) \subset L_2(0,1)$. Thus, an $X_\infty$-bounded function is $X_2$-bounded. The converse, however, is not true, as it is shown via an example in Ref. [17].

In this paper, it is assumed that the applied force $F(\cdot,\cdot)$ in Eq. (28) is $X_\infty$-bounded. The goal is to prove that the string displacement $y(\cdot,\cdot)$ is $X_\infty$-bounded; that is, to establish the BIBO stability of the string.

The stability result is as follows:

**Theorem 3.1:** Consider the string represented by Eqs. (28), (4), and (29). Let $F \in X_\infty$. The string displacement is bounded and satisfies

$$\|y\|_{X_\infty} \leq \left( \frac{K}{\varepsilon[(2\delta - \varepsilon K)(1 - v^2) - 2\delta v]} \right)^{1/2} \|F\|_{X_2} < \infty\,, \tag{30}$$

where

$$0 < \varepsilon < \frac{2\delta(1 - v - v^2)}{K(1 - v^2)}\,, \tag{31}$$

and $K$ is that in Eq. (11).

**Proof:** The proof can be established by taking steps similar to those is Ref. [17]; to avoid repetition, details are not presented. □

## 4. A Few Open Problems

In this section, a few open problems regarding the stability and stabilization of nonlinear moving strings with the Kelvin-Voigt damping are stated.

**Absence of viscous damping:** The stability results in Sections 2 and 3 crucially rely on the fact that there is viscous damping in the string, i.e., $\delta > 0$. Now, suppose that $\delta = 0$. In this case, the model of the moving string is

$$y_{tt}(x,t) + 2v\, y_{xt}(x,t) = \left[ 1 - v^2 + \frac{3}{2} b\, y_x^2(x,t) \right] y_{xx}(x,t)$$

$$+ 2\eta\, y_x(x,t)\, y_{xx}(x,t)\, y_{xt}(x,t) + \eta\, y_x^2(x,t)\, y_{xxt}(x,t)\,, \tag{32}$$

for all $x \in (0,1)$ and $t \geq 0$. For this string, let the boundary and initial conditions be, respectively, the same as those in Eqs. (4) and (5). In this string, the only damping is the Kelvin-Voigt damping $\eta > 0$.

Recall the function $E(\cdot)$ in Eq. (7). Clearly, this function is non-negative. Also recall that al least one of the functions $f$ or $g$ in Eq. (5) is not identically zero over $[0,1]$. Furthermore, the function $f$, for which $f(0) = 0$ by Eq. (4), cannot assume a non-zero constant value over $[0,1]$. Thus,

$$E(0) = \int_0^1 [\frac{1}{2}\, g^2(x) + \frac{1}{2}\,(1 - v^2)\, f_x^2(x) + \frac{b}{8}\, f_x^4(x)]\, dx > 0\,. \tag{33}$$

Computing the derivative of $t \mapsto E(t)$ with respect to $t$ along the solution of Eqs. (32), (4), and (5), and using Eqs. (14a)-(14d), it is concluded that

$$\dot{E}(t) = -\,\eta \int_0^1 y_x^2\, y_{xt}^2\, dx \leq 0\,, \tag{34}$$

for all $t \geq 0$. Since $\dot{E}(\cdot)$ is non-positive, the function $E(\cdot)$ is not increasing. However, the stability of the string cannot be established, unless LaSalle's invariance principle is used. It is known that using this principle for partial differential equations is prohibitively difficult; see, e.g., Refs. [18, 19]. Nevertheless, the following questions remain:

**Problem 4.1:** Is the moving string represented by Eqs. (32), (4), and (5), which has only the Kelvin-Voigt damping $\eta > 0$, stable?

**Problem 4.2:** Let an $X_\infty$-bounded distributed force be applied to the string represented by Eq. (32), and let the boundary and initial conditions, respectively, be the same as those in Eqs. (4) and (29). Is the forced moving string bounded-input bounded-output stable?

Affirmative answers to Problems 4.1 and 4.2 are desirable, because they prove that the Kelvin-Voigt damping is sufficient to guarantee the zero-input stability and the BIBO stability of the model.

**Boundary Control:** The stability of the string represented by Eqs. (3), (4), and (5) was established in Section 2 for fixed eyelets. It may be suggested that a boundary control at $x = 1$ would lead to faster convergence of the string displacement to zero or would relax condition in Eq. (20) for the stability. To apply the boundary control the eyelet on the right should be free to move transversally; see Fig. 2. The control input is proportional to $-y_t(1, t)$, namely, the negative feedback of the transversal velocity of the string at $x = 1$. In this case, the boundary conditions are

$$y(0, t) = 0, \qquad T(1, t)\, y_x(1, t) = -k_v\, y_t(1, t), \tag{35}$$

for all $t \geq 0$, where $T(1, \cdot)$ is the tension in the string at $x = 1$ and $k_v > 0$ is the control gain. A problem whose solution is desirable is as follows:

**Problem 4.3:** Is the string represented by Eqs. (3), (35), and (5) stable when (i) $\delta > 0$; (ii) $\delta = 0$? If yes, would a choice of $k_v$ result in faster convergence of the string displacement to zero?

## 5. Conclusions

In this paper, the stability of a nonlinear axially moving string with the Kelvin-Voigt damping was studied. First, it was proved that the transversal displacement of the string converges to zero when the axial speed of the string is less than a certain critical value. The proof was established by using an appropriate Lyapunov function, which tends to zero exponentially. Also, it was shown that the string is bounded-input bounded-output stable. In establishing these stability results, the viscous damping in the string model played an important role. Furthermore, a few open problems regarding the stability and stabilization of strings with the Kelvin-Voigt damping were stated. In these problems, the main quest is whether the Kelvin-Voigt damping can guarantee the zero-input stability and the BIBO stability of the string. Solutions of these problems would provide deeper understanding of the behavior of the nonlinear model of the string considered in this paper.

What was presented in this paper and further study of the model of the moving string with the Kelvin-Voigt damping would establish the validity of this model. Note that Eqs. (1), (2), and

(3) are only mathematical models of a real moving string. They would be reliable models if they predict the observed behavior of the real system. For instance, if the moving string is observed to be stable, then a proposed model must predict the stability. If the model fails to do so, then its validity is questionable.

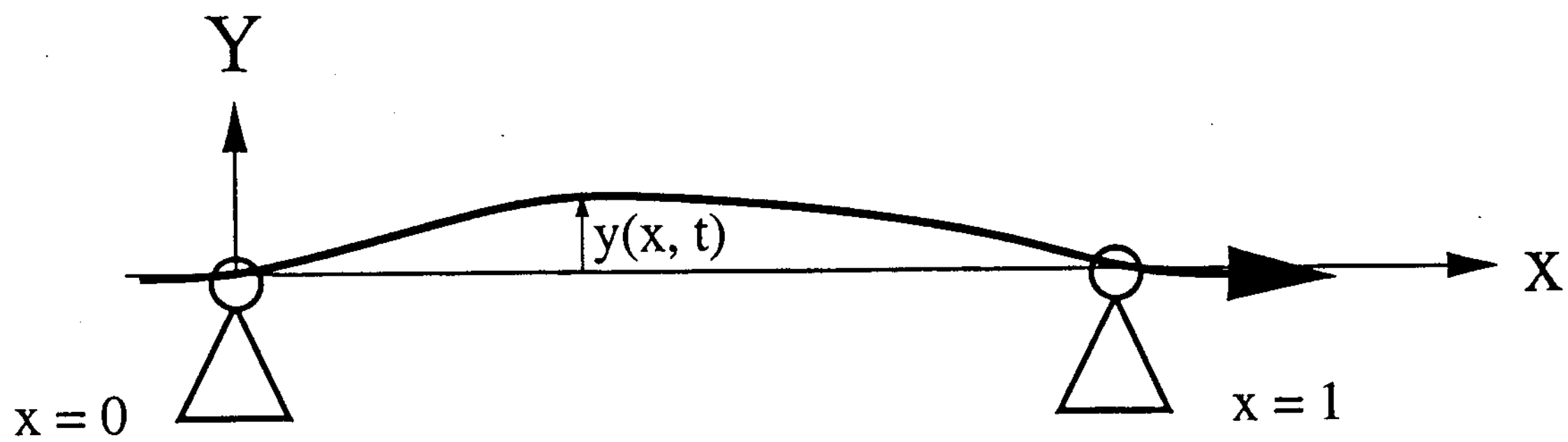

**Fig. 1.** A string-like continuum is pulled at a constant speed through two fixed eyelets. The distance between the eyelets is 1 .

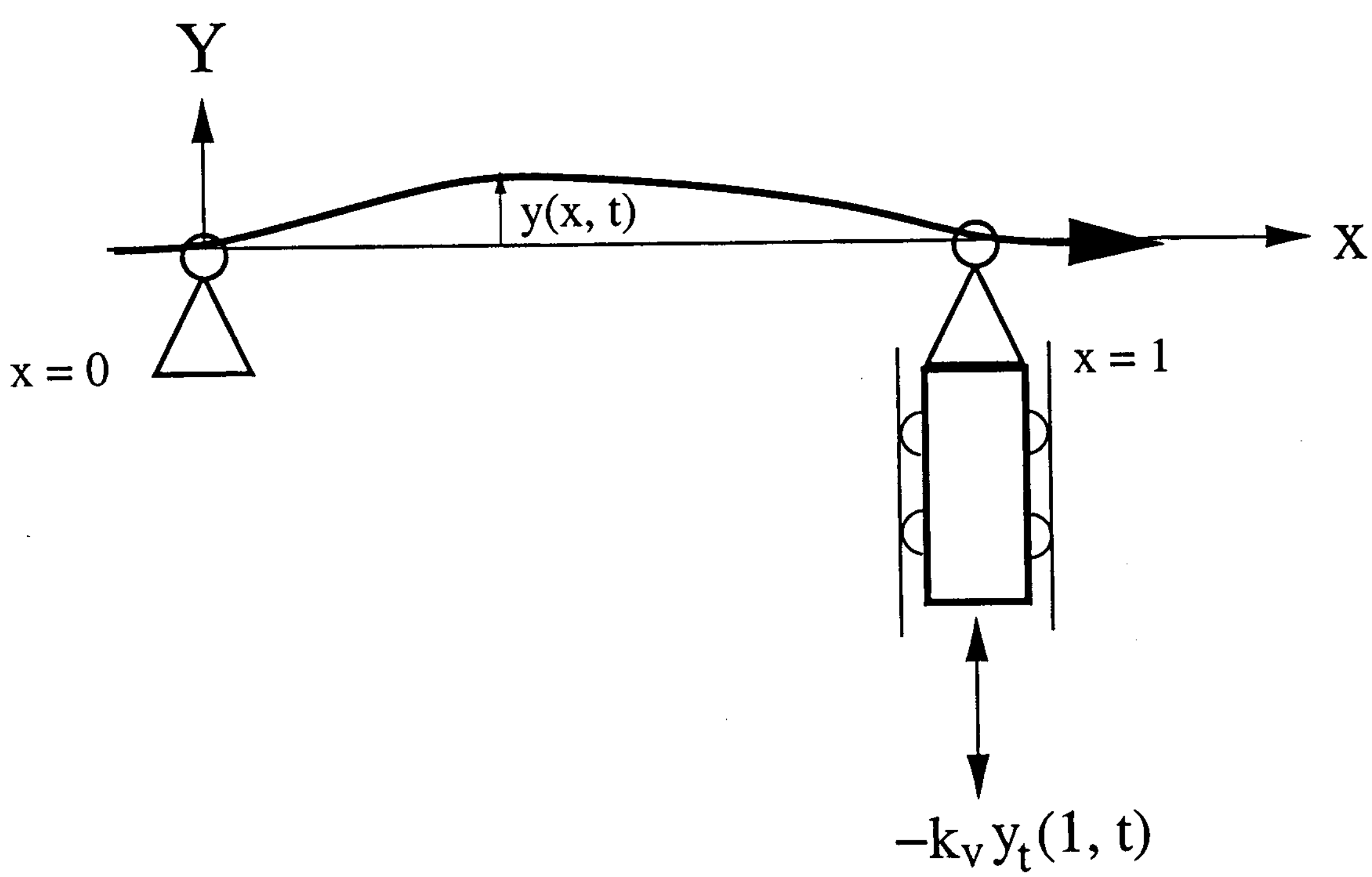

**Fig. 2.** The eyelet on the right is free to move transversally. A control input proportional to $-y_t(1, t)$ is applied to the string at $x = 1$ .